\documentclass[11pt]{amsart}

\usepackage{amssymb} 
\usepackage[mathscr]{eucal} 
\usepackage[matrix,arrow,tips,curve,ps]{xy}
\usepackage{amsmath} 
\usepackage{epsfig}
\usepackage{amscd}
\usepackage[usenames,dvipsnames,svgnames,table]{xcolor}
\usepackage{tikz}
\usepackage{tikz-cd}
\usepackage{mathrsfs}
\usepackage{verbatim}

\newtheorem{TEO}{Theorem}[section]

\newtheorem{REM}[TEO]{Remark}

\theoremstyle{definition}
\newtheorem{remark}[TEO]{Remark}
\newtheoremstyle{dico}
 {\baselineskip}   
  {\topsep}   
  {}  
  {0pt}       
  {} 
  {.}         
  {5pt plus 1pt minus 1pt} 
  {}          
\theoremstyle{dico}

\numberwithin{equation}{section}

\def\OO{{\mathcal O}}

\newcommand\dual{\mathrel{\raise3pt\hbox{$\underline{\mathrm{\thinspace d
\thinspace}}$}}}

\newcommand\proj{\mathbb P}
\newcommand\Z{\mathbb Z}

\newcommand\rank{\operatorname{rank}}

\newcommand{\restr}[1]          {\vert_{#1}}

\newcommand{\cinf}{C^\infty}

\renewcommand{\phi}{\varphi}

\newcommand{\cd}{\cdot}

\newcommand{\ra}{\rightarrow}

\newcommand{\debar }            {\bar {\partial } }

\newcommand{\mihi}[1]{}

\begin{document}

\footnote{
The first  author was partially supported by MIUR PRIN 2015 ``Geometry of Algebraic Varieties''.
   The second author was partially supported by MIUR PRIN
   2015 ``Moduli spaces and Lie theory'' and by FIRB 2012 `` Moduli Spaces and their Applications''.   The authors were also
   partially supported by GNSAGA of INdAM.
 }

\title{A bound on the dimension of a totally geodesic submanifold in the Prym locus}

\author[E. Colombo]{Elisabetta Colombo}
\address{Dipartimento di Matematica,
Universit\`a di Milano, via Saldini 50,
     I-20133, Milano, Italy } \email{{\tt
elisabetta.colombo@unimi.it}}

\author[P. Frediani]{Paola Frediani}
\address{ Dipartimento di Matematica, Universit\`a di Pavia,
via Ferrata 5, I-27100 Pavia, Italy } \email{{\tt
paola.frediani@unipv.it}}

\maketitle

\setlength{\parskip}{.1 in}

\begin{abstract}
We give an upper bound for the dimension of a germ of a totally geodesic submanifold, and hence of a Shimura variety of ${\mathcal A}_{g-1}$,  contained in the Prym locus. First we give such a bound for a germ passing through a Prym variety  of a $k$-gonal curve in terms of the gonality $k$. Then we deduce a bound only depending on the genus $g$.

\end{abstract}

\section{Introduction}
Denote by $ {\mathcal R}_g$ the moduli space parametrising isomorphism classes of pairs $[(C, \alpha)]$ where $C$ is a smooth complex projective curve of genus $g$ and $\alpha$ is a two-torsion  line bundle on $C$.
Consider the Prym map $Pr: {\mathcal R}_g \rightarrow {\mathcal A}_{g-1},$
which associates to a point $[(C,\alpha)] \in {\mathcal R}_g$ its Prym variety $P(C,\alpha)$ with its principal polarization.
We recall that the Prym map is generically an embedding for $g \geq 7$ (\cite{fs}, \cite{ka}).

Hence one can study the second fundamental form of the restriction of the Prym map to the open set ${\mathcal R}_g^0$ where the Prym map is an embedding, with respect to the orbifold metric  on ${\mathcal A}_{g-1}$ induced by the symmetric metric on the Siegel space $H_{g-1}$.  This has been done in \cite{cfprym} in analogy with the study of  the second fundamental form of the Torelli map $j: {\mathcal M}_g \rightarrow {\mathcal A}_{g}$ done in \cite{cpt} and \cite{cf2}.
In fact the expression of the second fundamental form both in the Torelli and the Prym case is very similar.

In the Torelli case, we recall that a conjecture by Coleman and Oort says that for  big enough genus,  there should not exist Shimura subvarieties of ${\mathcal A}_g$ generically contained in the Torelli locus, i.e. contained in $\overline {j({{\mathcal M}}_g)}$ and intersecting $j({\mathcal M}_g)$.

Since Shimura subvarieties are totally geodesic it makes sense to use the second fundamental form to investigate this problem. This has been done in \cite{cfg} where an upper bound on the dimension of a germ of a totally geodesic submanifold contained in ${\mathcal M}_g$, which only depends on $g$, is given.

In this paper, using the second fundamental form of the Prym map, we give a similar  upper bound on the dimension of a germ of a totally geodesic submanifold contained in the Prym locus depending only on $g$ (Theorem \ref{stima2}).

{\bf Theorem A} (See Theorem \ref{stima2})
  {\em If $g\geq 9$ and $Y$ is a germ of a totally geodesic submanifold of
  ${\mathcal A}_{g-1}$ contained in the Prym locus and intersecting $Pr({\mathcal R}_g^0)$, then $\dim Y \leq
  \frac{5}{2}g +\frac{1}{2}$.}

 This result follows from a bound depending on the gonality of the curve (Theorem \ref{stima1}).

{\bf Theorem B} (See Theorem \ref{stima1})
 {\em  Assume that $[(C,\alpha)] \in {\mathcal R_g^0}$ where $C$ is a $k$-gonal curve of genus $g$ with $g\geq 9$ and
  $k\geq 3$. Let $Y$ be a germ of a totally geodesic submanifold of
  ${\mathcal A}_{g-1}$ which is contained in the Prym locus and passes through
  $P(C, \alpha)$. Then $\dim Y \leq 2g+k - 1$.}

Motivated by these similarities in \cite{cfgp} we posed a similar question to the one of Coleman and Oort, namely we asked if, for big enough genus there exist Shimura subvarieties generically contained in the Prym locus. We say that a subvariety $Z \subset {\mathcal A}_{g-1}$ is generically contained in the Prym locus if $Z \subset \overline{Pr({\mathcal R}_g)}$, $Z \cap Pr({\mathcal R}_g) \neq \emptyset$ and $Z$ intersects the locus of irreducible principally polarised abelian varieties.

In \cite{cfgp} we explictitely constructed examples of Shimura curves  in
${\mathcal A}_{g-1}$ generically contained in the Prym locus, for $g \leq 13$.

\section{The 2nd fundamental form of the Prym map}

In this section we start by recalling the definition and some properties of the second fundamental form of the Prym map studied in \cite{cfprym}. Consider the Prym map $Pr: {\mathcal R}_g \rightarrow {\mathcal A}_{g-1},$
which associates to a point $[(C,\alpha)] \in {\mathcal R}_g$ its Prym variety $P(C,\alpha)$ with its principal polarization.
If $\pi: \tilde{C} \rightarrow C$ is the unramified double covering associated to the pair $(C,\alpha)$,
the Prym variety  of the double covering is the principally polarized abelian variety of dimension $g-1$
defined as the connected component containing the origin of the kernel of the norm map $Nm_{\pi}:J{\tilde C} \rightarrow JC$,
$$P(C,\alpha) = Ker(Nm_{\pi})^0 \subset J{\tilde C}.$$

On ${\mathcal A}_{g-1}$ there is an orbifold metric which is induced by the symmetric metric on the Siegel space $H_{g-1}$ of which ${\mathcal A} _{g-1}$ is a quotient by the action of $Sp(2g, \Z)$. We will consider ${\mathcal A}_{g-1}$ endowed with this orbifold metric.

We recall that the Prym map is generically an embedding for $g \geq 7$ (\cite{fs}, \cite{ka}). Hence there exists
an open set ${\mathcal R}^0_g \subset {\mathcal R_g}$ where $Pr$ is an embedding.

Consider the orbifold tangent bundle exact sequence of the Prym map
\begin{equation}
\label{tangent}
0 \rightarrow T_{{\mathcal R}^0_g} \rightarrow T_{{{\mathcal A}_{g-1}}_{|{\mathcal R}^0_g}} \rightarrow {\mathcal N_{{\mathcal R}^0_g/{\mathcal A}_{g-1}}} \rightarrow 0
\end{equation}

At the point $b_0:= [(C,\alpha)]$ the tangent bundle $T_{{{\mathcal A}_{g-1}}_{|{\mathcal R}^0_g}, b_0}$ is isomorphic to $S^2H^0(K_C \otimes \alpha)^*$,    $T_{{\mathcal R}^0_g,b_0}$ is isomorphic to
$ H^{1}(T_C)$ and the codifferential of the Prym map is the dual of the multiplication map. So the dual of the exact sequence \eqref{tangent} at the point $b_0$ becomes

$$0 \rightarrow I_2(K_C\otimes \alpha)  \rightarrow S^2H^0(K_C \otimes \alpha) \stackrel {m}\rightarrow H^0(K_C^{\otimes 2})
\rightarrow 0. $$

Therefore the dual of the second fundamental form of the Prym map at the point $b_0$  is a map

 \begin{equation}
 \label{II}
 II: I_2(K_C\otimes \alpha)\rightarrow H^0(2K_C)\otimes H^0(2K_C) \cong Hom(H^1(T_C), H^0(2K_C))
\end{equation}

In the proof of Theorem 2.1 of \cite{cfprym}  we showed that the second fundamental form $II$
is equal (up to scalar) to the Hodge-Gaussian map $\rho$ defined in Theorem 4.4 of \cite{cpt}. In particular $II$ is a lifting of the second Gaussian map of the line bundle  $K_C \otimes \alpha$ as
it happens for the second fundamental form of the period map and
the second Gaussian map of the canonical line bundle (see Theorems
2.1 and 4.5 of \cite{cpt}).
We will now recall the definition of the Hodge-Gaussian map $\rho$.

First of all recall  that given a
holomorphic line bundle $\alpha$ of degree zero on a curve $C$, there
exists a unique (up to constant) hermitian metric $H$ on $\alpha$ and a
unique connection $D_H$ on $\alpha$ which is compatible both with the
holomorphic structure and with the metric and which is flat. 

So we can write $D_H = D'_H +\overline{\partial}$, where $D'_H$ is the $(1,0)$ component. Such
a pair $(\alpha,H)$ is also called a harmonic line bundle which satisfies K\"ahler identities and a harmonic decomposition   (see \cite{sim}):
$${\mathcal A}^{\bullet}(\alpha) = {\mathcal H}^{\bullet}(\alpha) \oplus im(D_H) \oplus im(D^*_H)= {\mathcal H}^{\bullet}(\alpha) \oplus
im(\overline{\partial}) \oplus im(\overline{\partial}^*),$$ where
${\mathcal H}(\alpha)$ is the kernel of the laplacian operator $\Delta
= D_H D^*_H + D^*_H D_H =
2(\overline{\partial}\overline{\partial}^*+
\overline{\partial}^*\overline{\partial}),$
and the principle of two types
 $$ker(D'_H) \cap ker(\overline{\partial}) \cap (im(D'_H) + im(\overline{\partial})) = im(D'_H \overline{\partial}).$$

Let $v \in H^1(T_C)$ and $Q\in I_2(K_C \otimes \alpha)$. Take a representative $\theta_v \in A^{0,1}(C,T_C)$ of $v$, $[\theta_v ]= v \in H^1(T_C)$ and a basis $\{\omega_1,...,\omega_{g-1}\}$ of $H^0(K_C \otimes \alpha)$. Then $\theta_v(\omega_i) = \gamma_i + \overline{\partial}h_i$ where $\gamma_i \in A^{0,1}(C,\alpha)$ is harmonic and $h_i \in A^0(C,\alpha)$.

Let $Q = \sum_{i=1}^{g-1} a_{i,j} \omega_i \otimes \omega_j$ be a quadric in $I_2(K_C \otimes \alpha)$. The expression of $\rho$ given in Theorem 4.5 of  \cite{cpt} is the following

\begin{equation}
\label{rho}
 \rho(Q)(v)= \sum_{i,j} a_{i,j} \omega_j D'_H h_i,
 \end{equation}

 In Theorem 2.1 of  \cite{cfprym} we proved that $II(Q) = -2 \rho(Q)$.

In the case where the tangent vector $v$ is a Schiffer variation at a point $P \in C$ we have an equivalent expression of $\rho(Q)(v)$.
Schiffer variations can be seen intrinsically as the dual map $\xi: T_C^{2} \rightarrow H^1(T_C) \otimes {\mathcal O}_C$ of the evaluation map $H^0(K_C^2) \otimes {\mathcal O}_C \rightarrow K_C^2$.

Let us give a local description which is useful for the computations.

Consider a point $P\in C$, the coboundary of the exact sequence
  $0 \ra T_C \ra T_C(P) \ra T_C(P)\restr{P} \ra 0$ yields an injection
  $H^0(T_C(P)\restr{P}) \cong {\mathbb C}  \hookrightarrow H^1(C,T_C)$. A Schiffer variation at $P$  is an element of the image.  If
  $(U,z)$ is a chart centred at $P$ and $b\in \cinf_0(U)$ is a bump
  function which is equal to 1 on a neighbourhood of $P$, then
  \begin{gather*}
    \theta:= \frac{\debar b}{z}\cd \frac{\partial}{\partial z}
  \end{gather*}
  is a Dolbeault representative of a Schiffer variation $\xi_P$ at $P$.   It is well known
  that Schiffer variations generate $H^1(C,T_C)$ \cite[p.175]{acg}.

   Write locally $\omega_i = \phi_i(z) dz \otimes a$, where $a$ is a local section of $\alpha$. We have $\theta_P \cdot \omega_i = \overline{\partial}(\frac{b}{z} \phi_i a)= \gamma_i + \overline{\partial}h_i$ on $C - \{P\}$. So if we set $g_i(z) = \frac{b \phi_i}{z} a - h_i$ on $C-\{P\}$, we obtain that $\gamma_i = \overline{\partial}(g_i)$. Moreover in \cite{cpt} Theorem 4.5 it is proven that  the $(1,0)$-form on $C-\{P\}$ with values in $\alpha$, $\eta_i := D'_Hg_i$, satisfies $\overline{\partial} \eta_i = 0$ and

   $$\rho(Q) (\xi_P) = \sum_{i,j} a_{i,j}
\omega_i D'_H h_j =  -\sum_{i,j} a_{i,j} \omega_i\eta_j.$$

 In  \cite{cpt} Theorem 4.5 the authors also proved that all the $\eta_i$'s are proportional, namely, if $P$ is not a base point for $K_C \otimes \alpha$ and  $\phi_i(P) \neq 0$, then $\eta_j = \frac{\phi_j(P)}{\phi_i(P)} \eta_i$, $\forall j$.
In fact,
$\phi_i(P) \eta_j - \phi_j(P) \eta_i $ has no poles and since $\overline{\partial} \eta_i = 0$, it is holomorphic on $C$, hence it is harmonic. Moreover $\phi_i(P) \eta_j - \phi_j(P) \eta_i = D'_H(\phi_i(P)g_j - \phi_j(P) g_i),$ and $\phi_i(P)g_j - \phi_j(P) g_i \in A^0(C, \alpha)$.
So,  by the harmonic decomposition, we have $\phi_i(P) \eta_j - \phi_j(P) \eta_i=0$.

Finally, if we denote by $\eta_P : = \frac{\eta_i}{\phi_i(P)} $, we can compute  $\rho(Q) (\xi_P) = -\sum_{i,j} a_{i,j} \omega_i \phi_j(P) \eta_P$, so if we fix local coordinates and we evaluate at the point $R$,
we have
\begin{equation}
\label{P-R}
\xi_R(II(Q) (\xi_P))  = -4 \pi i (\rho(Q)(\xi_P))(R)= -  4 \pi i Q(P,R)  \cdot \eta_P(R) \ if  \ P \neq R
\end{equation}
 and

\begin{equation}
\label{P-P}
\xi_P(II(Q)(\xi_P)) = -4 \pi i (\rho(Q)(\xi_P))(P)=  -2 \pi i \mu_{2}(Q)(P),
\end{equation}
where $ \mu_{2}: I_2(K_C \otimes \alpha) \rightarrow H^0(4K_C)$ is the second Gaussian map of the bundle $K_C \otimes \alpha$. (cf. \cite{cfprym} Corollary 2.2).

\begin{REM}
\label{inj}
The second fundamental form $II: I_2(K_C \otimes \alpha) \rightarrow S^2H^0(2K_C)$ is injective.
\end{REM}
\proof
The proof is very similar to the proof of  {\cite[Corollary 3.4]{cf2}, namely, assume that $Q = \sum_{i=1}^{g-1} a_{i,j} \omega_i \otimes \omega_j$ is a quadric in $I_2(K_C \otimes \alpha)$ such that $II(Q) =0$. Then for very point $P \in C$ not in the base locus of $|K_C \otimes \alpha|$, we have $0= II(Q)(\xi_P) =  2\eta_P \sum_{i,j} a_{i,j} \omega_i \phi_j(P)$. Therefore $ \sum_{i,j} a_{i,j} \omega_i \phi_j(P) = 0$ and hence $\sum_j a_{i,j} \phi_j(P) =0$, $\forall i$ and for all $P \in C$ outside the base locus of $|K_C \otimes \alpha|$. Thus we must have $\sum_ja_{i,j} \omega_j =0$, $\forall i$  and so $a_{i,j} = 0$, $\forall i,j$.
\qed
 \section{Totally geodesic submanifolds}

In this section, following the ideas of \cite{cfg}, we give an upper bound for the dimension of a germ of a totally geodesic submanifold of $ {\mathcal A}_{g-1}$ contained in the Prym locus.

\begin{TEO}
  \label{rank}
  Assume that $[(C,A)] \in {\mathcal R_g^0}$ is such that $C$ is a  $k$-gonal curve of genus $g$, with $g\geq 9$
  and $k \geq 3$. Then there exists a quadric $Q \in I_2(K_C \otimes \alpha)$ such
  that $\operatorname{rank}II(Q) \geq  2g-2k -4$.
\end{TEO}
\proof Let us fix a local coordinate  at the relevant points and  write $\xi_P$ for a Schiffer
variation at $P$.  Let $F$ be a line bundle on $C$ such that $|F|$ is
a $g^1_k$ and choose a basis $\{x,y\}$ of $H^0(F)$.  Set $M = K_C \otimes \alpha \otimes F^{-1}$
and denote by $B$ the base locus of $|M|$, then by Clifford Theorem we have: $2(h^0(M(-B)) - 1) \leq deg(M(-B)) = 2g-2 - k - deg(B)$, hence $deg(B) \leq 2g-2 -k - 2h^0(M(-B)) +2 = 2g-k-2h^0(M) \leq k+2, $ by Riemann Roch.

Take a pencil $\langle t_1,t_{2}\rangle$ in $H^0(M)$, with base locus $B$.  Write $t_i = t'_i s$ for a
section $s\in H^0(C,\OO_C(B)) $ with $\operatorname{div} (s) = B$.
Then $\langle t'_1, t'_2\rangle $ is a base point free pencil in $|M(-B)|$.  Let
$\psi : C \rightarrow  \proj^1$ be the morphism induced by this pencil.

Consider the rank 4 quadric $Q:= xt_1 \odot yt_2 - xt_2 \odot
yt_1$. Clearly $Q \in I_2(K_C \otimes \alpha)$. Set $d:= \deg(M(-B)) = 2g-2-k-
\deg(B)$.  As in the proof of \cite{cfg} Theorem 4.1 one can  show that if
  $\{P_1,...,P_d\}$ is a fibre of the morphism $\psi$ over a regular
value, then the Schiffer variations $\xi_{P_1},...,
\xi_{P_d}$ are linearly independent in $H^1(C,T_C)$.  Denote by $\phi$ the morphism induced
by the pencil $|F|$ and consider the set $E:= \psi(Crit(\phi) \cup
Crit(\psi) \cup B)$ where $Crit(\phi) $ (resp.  $Crit(\psi) $)
denote the set of critical points of $\phi$ (resp. $\psi$).  Let $z
\in {\proj}^1 \setminus E$ and let
$\{P_1,\ldots,P_d\}$ be the fibre of
$\psi$ over $z$. By changing coordinates on $\proj^1$ we can assume $z =
[0,1]$, i.e. $t'_1(P_i) = 0$ for $i =1,\ldots d$. Then clearly
$t_1(P_i) = 0$, so $Q(P_i,P_j) = 0$ for all $i,j$.  By \eqref{P-R},\eqref{P-P}  one immediately obtains that the restriction of
$II(Q)$ to the subspace $W:= \langle \xi_{P_1},...,\xi_{P_d} \rangle$ is represented in the basis
$\{\xi_{P_1},..., \xi_{P_d}\}$ by a diagonal matrix with entries $-2 \pi
i \mu_2(Q)(P_i)$ on the diagonal.  For a rank 4 quadric the second
Gaussian map can be computed as follows: $ \mu_2(Q) = \mu_{1,F}(x
\wedge y)\mu_{1,M}(t_1 \wedge t_2)$, where $\mu_{1,F}$ and $\mu_{1,M}$ are the first Gaussian maps of the line bundles $F$ and $M$ (see \cite[Lemma 2.2]{cf1}).  Now
$\mu_{1,F}(x\wedge y) (P_i) \neq 0$, because $P_i \not \in
Crit(\phi)$ by the choice of $z$.  Moreover
$P_i \not \in B$.  On $C\setminus B$ the morphism $\psi$ coincides
with the map associated to $\langle t_1, t_2\rangle$. Since $P_i \not \in
Crit(\psi)$, it is not a critical point for the latter map. Therefore
also $\mu_{1,M}(t_1\wedge t_2) (P_i) \neq 0$.
Thus $\mu_2(Q)(P_i) = \mu_{1,F}(x \wedge y)(P_i)\mu_{1,M}(t_1 \wedge
t_2)(P_i) \neq 0$ for every $i=1,..., d$.  This shows that in the
basis $\{\xi_{P_1},..., \xi_{P_d}\}$ the quadric $II(Q)_{|W}$ is
represented by a diagonal matrix with non-zero diagonal entries.  So
$II(Q)$ has rank at least $d = 2g-2-k-\deg B \geq  2g-2k-4$.\qed

\begin{TEO}
  \label{stima1}
  Assume that $[(C,\alpha)] \in {\mathcal R_g^0}$ where $C$ is a $k$-gonal curve of genus $g$ with $g\geq 9$ and
  $k\geq 3$. Let $Y$ be a germ of a totally geodesic submanifold of
  ${\mathcal A}_{g-1}$ which is contained in the Prym locus and passes through
  $P(C, \alpha)$. Then $\dim Y \leq 2g+k - 1$.
\end{TEO}
\proof By Theorem \ref{rank} we know that there exists a quadric $Q
\in I_2(L)$ such that the rank of $II(Q)$ is at least $2g-2k-4$. By
assumption for any $v \in T_{[C]}Y$ we must have that $v(II(Q)(v)) = II(Q)(v \odot v) = 0$, so $v$ is isotropic for $II(Q)$, hence
\begin{gather*}
  \dim T_{[C]} Y \leq 3g-3 - \frac{(2g-2k-4) }{2}= 2g+ k - 1
  .
\end{gather*}
\qed
\begin{remark}
  In Theorem \ref{rank} if $|M|$ is base point free the
  above proof shows that $ \rank II(Q) \geq 2g-2-k$.  So in this case the bound on the dimension of  a
  germ of a totally geodesic submanifold  $Y$  of ${\mathcal A}_{g-1}$ contained in the
  Prym locus and passing through $P(C, \alpha)$ with $C$ a $k$-gonal curve such that $|M|$ is base point free becomes:
  $\dim Y \leq 2g-2 + {k}/{2}$.
\end{remark}

\begin{TEO}
  \label{stima2}
  If $g\geq 9$ and $Y$ is a germ of a totally geodesic submanifold of
  ${\mathcal A}_{g-1}$ contained in the Prym locus and intersecting $Pr({\mathcal R}_g^0)$, then $\dim Y \leq
  \frac{5}{2}g +\frac{1}{2}$.
\end{TEO}
\proof This immediately follows from Theorem \ref{stima1}, since the
gonality of a genus $g$ curve is at most $ [(g+3)/{2}]$.  \qed

\begin{remark}
  As it is evident from the proof, gonality is used to construct a
  quadric $Q \in I_2(K_C \otimes \alpha )$ of rank 4 such that $II(Q)$ has large
  rank. It seems unlikely that gonality plays any role in this
  problem.  In fact we expect the existence of $Q \in I_2(K_C \otimes \alpha)$ with
  image $II(Q)$ a nondegenerate quadric on $H^1(C,T_C)$. This would
  give the upper bound $\frac{3}{2}(g-1)$ for the dimension of a germ
  of a totally geodesic submanifold.  On the other hand, the
    map $II$ is injective (see \ref{inj}). Therefore $II(I_2(K_C \otimes \alpha))$ gives a linear
    system of quadrics of dimension $\frac{(g-1)(g-6)}{2}$ on
    $\proj(H^1(C,T_C)) =\proj^{3g-4}$.  This already gives an upper bound
    for the dimension of a submanifold $Y$ as in Theorem
    \ref{stima2}. Indeed for any point $[(C,\alpha)]\in Y$, the tangent space
    $T_{[(C,\alpha)]} Y \subset H^1(C,T_C)$ is contained in the base locus of
    $II(I_2(K_C \otimes \alpha))$.  This means that $II(I_2(K_C \otimes \alpha))$ is contained
    in the space of quadrics $ q \in S^2H^0(C,2K_C)$ that vanish on
    $T_{[(C,\alpha)]}Y$. This yields the bound
    \begin{gather*}
      \dim Y \leq \frac{ -1 + \sqrt{ 32g^2 -32g +1}} {2}.
    \end{gather*}
    Nevertheless for any $g\geq 9$ this bound is weaker than the one
    provided in Theorem \ref{stima2}.  At any case the study of the
    base locus of the linear system $II(I_2(K_C \otimes \alpha))$ should clearly
    improve the understanding of totally geodesic submanifolds.  If
    one could prove that the base locus is empty, one would rule out
    the existence of totally geodesic submanifolds passing through
    $P(C,\alpha)$.

\end{remark}

\begin{remark}
Consider a family $(C_t, \alpha_t)$ contained in ${\mathcal R}^0_g$ and denote by $Z$ its image in ${\mathcal A}_{g-1}$ via the Prym map. Assume that there exists a group $G$ acting on the curves $C_t$ and admitting a linearization of the action on the  line bundles $\alpha_t$.
Then, arguing as in Lemma 5.1 of \cite{cfg}, one can easily prove that  the second fundamental form of $Z$ in  ${\mathcal A}_{g-1}$  is $G$-equivariant.

From this follows as in  Proposition 5.2 of \cite{cfg}  that if $I_2(K_{C_t} \otimes \alpha_t)^G =0$, then $Z$ is a totally geodesic subvariety of ${\mathcal A}_{g-1}$ contained in the Prym locus.

The examples constructed in  \cite{cfgp} that are contained in ${\mathcal R}_g^0$ satisfy the above property and are totally geodesic.
\end{remark}

\end{document}